\documentclass[a4paper,12pt]{article}
\usepackage[all]{xy}
\usepackage{amsmath, amsthm}
\usepackage{amssymb,amsfonts}
\usepackage{array}
\usepackage{hyperref}

\title{Characteristic polynomial and Wiener index of the Compressed Zero Divisor graph $\Gamma_{E}[\mathbb {Z}_{p^n}]$}
\author{B.Surendranath Reddy,Rupali.S.Jain and N.Laxmikanth\\
\textit{surendra.phd@gmail.com},{rupalisjain@gmail.com} and\\
{laxmikanth.nandala@gmail.com} }{}
\date{}

\begin{document}
\theoremstyle{definition}
\newtheorem{definition}{Definition}[section]
\newtheorem{example}[definition]{Example}
\newtheorem{remark}[definition]{Remark}
\newtheorem{observation}[definition]{Observation}
\theoremstyle{plain}
\newtheorem{theorem}[definition]{Theorem}
\newtheorem{lemma}[definition]{Lemma}
\newtheorem{proposition}[definition]{Proposition}
\newtheorem{corollary}[definition]{Corollary}
\newtheorem{AMS}[definition]{AMS}
\newtheorem{keyword}[definition]{keyword}
\maketitle
\section*{Abstract}
The Zero divisor Graph of a commutative ring $R$, denoted by $\Gamma[R]$, is a graph whose vertices are non-zero zero divisors of $R$ and two vertices are adjacent if their product is zero. In this paper we derive the Wiener index and the characteristic polynomial  of the Compressed zero divisor graph $\Gamma_{E}[\mathbb{Z}_m]$ where $m=p^n$ with prime $p$.\\ 
Keywords:-Compressed zero divisor graph,Wiener index and  characteristic polynomial 
\section{Introduction}
In this article,  section 2, is about the preliminaries and notations related to zero divisor graph of a commutative ring $R$ and compressed zero divisor graph, in section 3, we derive the characteristic polynomial of the compressed zero divisor  graph  $\Gamma_{E}[{\mathbb{Z}_{p^n}}]$, and in section 4, we calculate the  Wiener index of the compressed zero divisor  $\Gamma_{E}[{\mathbb{Z}_{p^n}}]$.
 \section{Preliminaries and Notations}
\begin{definition}{\textbf{Zero divisor Graph}}{\cite{1,2}}\\
	Let R be a commutative ring with unity and $ Z[R]$ be the set of its zero divisors. Then the zero divisor graph of R denoted by  $\Gamma[R]$, is the graph(undirected) with vertex set  $Z^*[R]= Z[R]-\{\mathbf{0}\}$, the non-zero zero divisors of $R$, such that two vertices $v,w \in Z^*[R]$ are adjacent if $vw=0$.
\end{definition}
\begin{definition}{\textbf{Adjacency matrix of $\Gamma[R]$}}{\cite{4}}\\
	The adjacency matrix of the zero divisor graph $\Gamma[R]$ is the matrix $[v_{ij}]$ with rows and columns labeled by the vertices and is given by
	\begin{align*} v_{ij}=\left\{
		\begin{array}{ll}
			1, & \hbox{$v_i\, \text{is adjacent to}\, v_j$;} \\
			0, & \hbox{$otherwise$.}
		\end{array} \right.
	\end{align*}
	Adjacency matrix of $\Gamma[R]$ is denoted by $M(\Gamma[R])$. Clearly for an undirected graph, the adjacency matrix is symmetric.
\end{definition}
\begin{definition}{\textbf{Compressed zero divisor graph $\Gamma[R]$}}{\cite{6,8}}\\
Let $R$ be a commutative ring with non-zero identity.Then the compressed zero divisor graph $\Gamma_{E}[R]$ of $R$ is defined by using the equivalence relation $\sim$ on $R$ given by x$\sim$y iff $ann_{R}(x) = ann_{R}(y)$.\\The set of vertices of $\Gamma_{E}[R]$ is V($\Gamma_{E}[R]$) =$\{[x]_{\sim} /x \in R / [0]_{\sim}\cup [1]_{\sim}\}$  and two distinct vertices $[x]_{\sim}$ and $[y]_{\sim}$ are connected by an edge in $\Gamma_{E}[R]$ iff $xy =0$ 
\end{definition}
\begin{definition}{\textbf{Wiener index of $\Gamma[R]$}}\\
	Let  $\Gamma[R]$ be a zero divisor graph with vertex set $V$. We denote the length of the shortest path between every pair of vertices $x,y\in V$ with d(x,y). Then the Wiener index of  $\Gamma[R]$ is the sum of the distances between all pair of vertices of  $\Gamma[R]$, i.e.,  $ W(\Gamma[R])= \sum\limits_{x,y\in V} d(x,y) $.\\
\end{definition}
\section{Characteristic polynomial of $\Gamma_{E}[\mathbb{Z}_{p^n}]$}
For any elements x and y of $R$, we define x$\sim$y iff $ann_{R}(x) = ann_{R}(y)$ where $\sim$ equivalence relation  on $R$. \\ The operation on the equivalence classes is defined as $[x]_{R}.[y]_{R}$=  $[xy]_{R}$ is well-defined.\\
 $\Gamma_{E}[R]$ denote the compressed zero divisor graph whose vertices are the elements from commutative ring with non-zero identity other than $[0]$ and  $[1]$ and two distinct vertices are adjacent iff  $[x]_{R}.[y]_{R}$=  $[0]_{R}$ provided $xy=0 $ note that x and y are distinct vertices in  $\Gamma[R]$ , so that $[x]_{R}$,$[y]_{R}$ are adjacent in  $\Gamma_{E}[R]$ iff $[x]_{R}$$\neq$ $[y]_{R}$.\\
In this section we derive the standard form of adjacency matrix of $\Gamma_{E}[\mathbb{Z}_{p^{n}}]$  and we find the corresponding characteristic polynomials.\\To start with we consider $n=6$.\\

Let  $n=6$
i.e., $\Gamma_{E}[\mathbb{Z}_{p^6}]$ the vertex set is given by\\
V( $\Gamma_{E}[\mathbb{Z}_{p^6}]$) = $\{a=m[p],b=m[p^2],c=m[p^3],d=m[p^4],e=m[p^5]\} $
and the adjacency matrix is given by 
$M(\Gamma[{\mathbb{Z}_{p^6}}])=\begin{bmatrix}
	0 & 0 & 0 & 0 & 1\\
	0 & 0 & 0 & 1 & 1\\
	0 & 0 & 1 & 1 & 1\\
	0 & 1 & 1 & 1 & 1\\
	1 & 1 & 1 & 1 & 1
\end{bmatrix} $\\
Then  $|M-\lambda{I}| = \begin{vmatrix}
	-\lambda & 0 & 0 & 0 & 1 \\
	0 & -\lambda & 0 & 1 & 1 \\
	0 & 0 & 1-\lambda & 1 & 1 \\
	0 & 1 & 1 & 1-\lambda & 1 \\
   1 & 1 & 1 & 1 & 1-\lambda  
	
\end{vmatrix}$\\
on expanding the determinant we get the characteristic polynomial as\\ 
$\lambda^5-3\lambda^4-3\lambda^3+4\lambda^2+\lambda-1$.\\
Let  $n=7$
i.e., $\Gamma_{E}[\mathbb{Z}_{p^7}]$ the vertex set is given by\\
V( $\Gamma_{E}[\mathbb{Z}_{p^7}]$) = $\{a=m[p],b=m[p^2],c=m[p^3],d=m[p^4],e=m[p^5],f=m[p^6]\} $
and the adjacency matrix is given by 
$M(\Gamma[{\mathbb{Z}_{p^7}}])=\begin{bmatrix}
	0 & 0 & 0 & 0 & 0 & 1\\
	0 & 0 & 0 & 0 & 1 & 1\\
	0 & 0 & 0 & 1 & 1 & 1\\
	0 & 0 & 1 & 1 & 1 & 1\\
	0 & 1 & 1 & 1 & 1 & 1\\
	1 & 1 & 1 & 1 & 1 & 1
\end{bmatrix} $\\
Then  $|M-\lambda{I}| = \begin{vmatrix}
	-\lambda & 0 & 0 & 0 & 0 & 1 \\
	0 & -\lambda & 0 & 0 & 1 & 1 \\
	0 & 0 & -\lambda & 1 & 1 & 1 \\
	0 & 0 & 1 & 1-\lambda & 1 & 1 \\
	0 & 1 & 1 & 1 & 1-\lambda & 1 \\
	1 & 1 & 1 & 1 & 1 & 1-\lambda  
	
\end{vmatrix}$\\
on expanding the determinant we get the characteristic polynomial as\\ 
$\lambda^6-3\lambda^5-6\lambda^4+4\lambda^3+5\lambda^2-\lambda+1$.\\
From the above examples and many more, we have observed that the coeffients of the characteristic polynomial form a pascal like traingle which is shown below\\
\[\begin{array}{ccccccccccccccccc}
&&& &    &    &    &    &  1 &    &    &    &    &  &&& \cr
&&& &    &    &    &  1 &    &  1 &    &    &    &  && &\cr
&& &&    &    &  1 &    &  1 &    &  1 &    &    &  & &&\cr
& &&&    &  1 &    &  2 &    &  1 &    &  1 &    &  & &&\cr
&& & & 1 &    &  2 &    &  3 &    &  1 &    &  1 &  & &&\cr
&& &1 &    &  3 &    & 3 &    & 4 &    &  1 &    & 1 &&&\cr
&&1&    &  3  &   &  6  &  &  4  &  &   5 &   &  1  &   &1&&\cr
&1&& 4   &    & 6  &    & 10 &    & 5 &    & 6  &    &  1 &&1&\cr
1&& 4&    &  10  &   &  10  &  & 15   &  & 6   &   & 7   &   &1&&1\cr
 \dots \cr
 \end{array}\]
This motivate us to derive the characteristic polynomial for $\Gamma[{\mathbb{Z}_{m}}]$, where $m=p^n$ with $p$ prime.

\begin{theorem}
Let $m=p^n$ with $p$ prime. Then the characteristic equation of  a compressed zero divisor graph  $\Gamma[{\mathbb{Z}_{m}}]$ is \\
$\lambda^{n-1}-b_{1}{\lambda^{n-2}}-b_{2}{\lambda^{n-3}} +b_{3}{\lambda^{n-4}} +........-(-1)^{\lfloor\frac{k+1}{2}\rfloor}[b_{k}{\lambda^{n-(k+1)}}]+.........+(-1)^{b_{1}}[b_{n-1}]\lambda^{n-n} = 0$.\\
where $b_{i}$= $p_{i}\choose{i} $ and $p_{i}=  \big{\lfloor\frac{n-1+i}{2}\rfloor} =floor function $, for $i= 1,2,3,....n-1$.
\begin{proof}
Let $m=p^n$. Then the set of non-zero zero divisors of $\mathbb{Z}_m$ is\\
$Z^*[\mathbb{Z_\mathbf{m}}]=\{[p],[p^2],[p^3],....[p^{n-1}]\}$ with cardinality  $n-1$
and the adjacency matrix is given by \\
$M(\Gamma[{\mathbb{Z}_{p^n}}])=\begin{bmatrix}
0 & 0 & 0 & 0 &\cdots& 0 & 1\\
0 & 0 & 0 & 0 & \cdots&1 & 1\\
\vdots& \vdots& \vdots&\vdots&\cdots&\vdots&\vdots\\
0 & 0 & 1&  1&\cdots & 1 & 1\\
0 & 1 & 1 & 1&\cdots & 1 & 1\\
1 & 1 & 1 & 1&\cdots & 1 & 1
\end{bmatrix} $\\

$|\mathbf{M}-\lambda{I}| =\begin{vmatrix}
-\lambda & 0 & 0 & 0 &\cdots& 0 & 1\\
0 & -\lambda & 0 & 0 & \cdots&1 & 1\\
\vdots& \vdots& \vdots&\vdots&\cdots&\vdots&\vdots\\
0 & 0 & 1&  1&\cdots & 1 & 1\\
0 & 1 & 1 & 1&\cdots & 1-\lambda & 1\\
1 & 1 & 1 & 1&\cdots & 1 & 1-\lambda
\end{vmatrix} $\\ 

Then $|\mathbf{M}-\lambda{I}|=|A||D-CA^{-1}B|=0$.\\
Since  $A$ is a null matrix of order $\big{\lfloor\frac{n-1}{2}\rfloor}$, we get  $|A|=(-\lambda)^{\big{\lfloor\frac{n-1}{2}\rfloor}}$.\\
And  $ |CA^{-1}B|=
\begin{vmatrix}
0 & 0 & 0 &\cdots & 0 \\
0 & \frac{-1}{\lambda} &\frac{-1}{\lambda} & \cdots & \frac{-1}{\lambda} \\
0 & \frac{-1}{\lambda} &\frac{-2}{\lambda} & \cdots & \frac{-2}{\lambda} \\
\vdots  & \vdots &\vdots & \ddots & \vdots  \\
0 & \frac{-1}{\lambda}& \frac{-2}{\lambda}&\cdots & \frac{{-\big{\lfloor\frac{n-1}{2}\rfloor}}}{\lambda}
\end{vmatrix}$\\[2mm]
 So, in general  when n is odd $ |D-CA^{-1}B|$ is given by\\
  $ |D-CA^{-1}B|=\begin{vmatrix}
 a_{3}-\lambda & a_{3} & a_{3}&\cdots & a_{3} \\
 a_{3} &  a_{5} -\lambda &  a_{5} &\cdots &  a_{5}  \\
  a_{3} & a_{5}& a_{7}-\lambda&\cdots & a_{7}\\   
 \vdots  & \vdots  & \ddots & \vdots  \\
  a_{3}  &  a_{5} & a_{7} & \cdots &  a_{n} -\lambda
 \end{vmatrix}$\\[2mm]
 
  when n is even $ |D-CA^{-1}B|$ is given by\\
 $ |D-CA^{-1}B|=\begin{vmatrix}
 a_{2}-\lambda & a_{2} & a_{2}&\cdots & a_{2} \\
 a_{2} &  a_{4} -\lambda &  a_{4} &\cdots &  a_{4}  \\
 a_{2} & a_{4}& a_{6}-\lambda&\cdots & a_{6}\\   
 \vdots  & \vdots  & \ddots & \vdots  \\
 a_{2}  &  a_{4} & a_{6} & \cdots &  a_{n} -\lambda
 \end{vmatrix}$\\[2mm]
 where $a_{n}=1+\frac{\big{\lfloor\frac{n-1}{2}\rfloor}}{\lambda}$ in both the cases for $n=2,3,4.... $.\\
 on solving the above determinant and calculating $|A||D-CA^{-1}B|$\\
 The characteristic equation of  a compressed zero divisor graph  $\Gamma[{\mathbb{Z}_{m}}]$ is \\
 $\lambda^{n-1}-b_{1}{\lambda^{n-2}}-b_{2}{\lambda^{n-3}} +b_{3}{\lambda^{n-4}} +........-(-1)^{[\frac{k+1}{2}]}[b_{k}{\lambda^{n-(k+1)}}]+.........+(-1)^{b_{1}}[b_{n-1}]\lambda^{n-n} = 0$.\\
 where $b_{i}$ =  ${p_{i}}\choose {i} $ and $p_{i}=  \big{\lfloor\frac{n-1+i}{2}\rfloor}$, for $i= 1,2,3,....n-1$.
\end{proof}
\end{theorem}

\section{ Wiener index of a Compressed zero divisor graph  $\Gamma_{E}[{\mathbb{Z}_n}]$. }
In this section, we derive the standard form of the Wiener  index of the compressed zero divisor graph $\Gamma_{E}[\mathbb{Z}_m]$ for $m=p^n$.  
To start with, we consider two cases $n$ being even and odd. 
\begin{theorem}\label{r1}
	Let $m=p^6$ with $p$ prime. Then the Wiener index of  a compressed zero divisor graph  $\Gamma_{E}[{\mathbb{Z}_{m}}]$ is 
	$\mathbf W(\Gamma_{E}[{\mathbb{Z}_{m}}])=14.$
	\begin{proof}
	Let $m=p^6$. Then the set of non-zero zero divisors of $\mathbb{Z}_m$ is\\
	$Z^*[\mathbb{Z_\mathbf{m}}]=\{[p],[p^2],[p^3],[p^4],[p^5]\}$ with cardinality  $n-1$. \\
	Let $ a=[p],b=[p^2] , c=[p^3],d=[p^4],e=[p^5].$\\
	Since the element a is not adjacent with b,c and d \\
	therefore $ d(a,b)=d(a,c)=d(a,d)=2 $ and $ d(a,e)=1 $\\
	and b is also not adjacent with  c therefore $ d(b,c)=2 $\\
	c is adjacent with d and e therefore $ d(c,d)=1$ and  $ d(c,e)= 1$\\
	but d adjacent with b,c e therefore $ d(d,b)=1 $ ,$ d(d,e)=1 $
		also e is adjacent with all other elements  except a therefore $ d(e,b)=1 $ \\
	$\therefore  \mathbf W(\Gamma[{\mathbb{Z}_{p^{4}}}])=\sum d(x,y) =d(a,b)+d(a,c)+d(a,d)+d(a,e)+d(b,c)+d(b,d)+d(b,e)+d(c,d)+d(c,e)+d(d,e)=2+2+2+1+2+1+1+1+1+1=14.$\\	
	\end{proof}
\end{theorem}
\begin{theorem}\label{r1}
	Let $m=p^7$ with $p$ prime. Then the Wiener index of  a compressed zero divisor graph  $\Gamma_{E}[{\mathbb{Z}_{m}}]$ is 
	$\mathbf W(\Gamma_{E}[{\mathbb{Z}_{m}}])=21.
	$.
	\begin{proof}
		Let $m=p^7$. Then the set of non-zero zero divisors of $\mathbb{Z}_m$ is\\
	$Z^*[\mathbb{Z_\mathbf{m}}]=\{[p],[p^2],[p^3],[p^4],[p^5],[p^6]\}$ with cardinality  $n-1$. \\
	Let $ a=[p],b=[p^2] , c=[p^3],d=[p^4],e=[p^5],f=[p^6].$\\
	Since the element a is not adjacent with b,c,d and e and it is only adjacent with f \\
	therefore $ d(a,b)=d(a,c)=d(a,d)=d(a,e)=2  $ and $ d(a,f)=1 $\\
	and b is also not adjacent with  c,d but  adjacent with e, f therefore $ d(b,c)=d(b,d)=2$ and $ d(b,e)=d(b,f)=1$\\
	c is adjacent with d,e and f therefore $ d(c,d)=1$,$ d(c,e)=1$ and  $ d(c,f)= 1$\\
	but d adjacent with e,f therefore $ d(d,e)=1 $ ,$ d(d,f)=1 $
	also e is adjacent with all other elements  except a therefore $ d(e,f)=1 $ \\
	$\therefore  \mathbf W(\Gamma_{E}[{\mathbb{Z}_{p^{7}}}])=\sum d(x,y) =d(a,b)+d(a,c)+d(a,d)+d(a,e)+d(a,f)+d(b,c)+d(b,d)+d(b,e)+d(b,f)+d(c,d)+d(c,e)+d(c,f)+d(d,e)+d(d,f)+d(e,f)=2+2+2+2+1+2+2+1+1+1+1+1+1+1+1=21.$\\	    	
	\end{proof}
\end{theorem}
With the above proofs we now head towards a general case
\begin{theorem}\label{r1}
Let $m=p^n$ with $p$ prime. Then the Wiener index of  a compressed zero divisor graph  $\Gamma_{E}{\mathbb{Z}_{m}}]$ is 
$\mathbf W(\Gamma_{E}[{\mathbb{Z}_{m}}])$=
$\frac{(n-2)(3n-4)}{2}$ if n is even and
$\mathbf W(\Gamma_{E}[{\mathbb{Z}_{m}}])$=
$\frac{(n-1)(3n-7)}{2}$ if n is odd.
 \begin{proof}
Let $m=p^n$. Then the set of non-zero zero divisors of $\mathbb{Z}_m$ is\\
$Z^*[\mathbb{Z_\mathbf{m}}]=\{[p],[p^2],[p^3],....[p^{n-1}]\}$ with cardinality  $n-1$. \\
Let $ a_{1}=[p],a_{2}=[p^2] , a_{3}=[p^3],.......a_{n-1}=[p^{n-1}].$\\
Since the element $a_{1}$ is not adjacent with $a_{2}$,$a_{3}$,....$a_{n-2}$ and it is only adjacent with $a_{n-1}$ \\
therefore $ d(a_{1},a_{2})=d(a_{1},a_{3})=......d(a_{1},a_{n-2})=2 $ and $d(a_{1},a_{n-1})= 1$\\
with a similar argument, we get\\
  $ d(a_{2},a_{3})=d(a_{2},a_{4})=......d(a_{2},a_{n-3})=2 $\\ and $d(a_{2},a_{n-2})= 1,d(a_{2},a_{n-1})= 1$\\
$\therefore $  $ d(a_{i},a_{i+1})=d(a_{i},a_{i+2})=......d(a_{i},a_{n-i-1})=2 $ and \\$d(a_{i},a_{n-i})= 1,d(a_{i},a_{n-i+1})= ......d(a_{i},a_{n-1})=1$\\
 so on $d(a_{n-1},a_{n})=1$\\
The combinations of the elements with distance is 2 are\\
$ d(a_{1},a_{2})=d(a_{1},a_{3})=......d(a_{1},a_{n-2})=2 $ \\
these are (n-3)combinations as $a_{1}$ and  $a_{n-1}$ are to be neglected.\\
similarly $ d(a_{2},a_{3})=d(a_{2},a_{4})=......d(a_{2},a_{n-3})=2 $ \\
these are (n-5)combinations as $a_{1}$, $a_{2}$ , $a_{n-2}$  and  $a_{n-1}$ are to be neglected.\\
so on $d(a_{ \lfloor \frac{n}{2}\rfloor-1},a_{\lfloor\frac{n}{2}\rfloor})=2$ this is only one combination as the suffixes are consecutive terms.\\ where $\lfloor \frac{n}{2}\rfloor= floor function$.\\
Now the combinations of the elements with distance is 1 are\\
$d(a_{1},a_{n-1})= 1$ this is only one term,\\
$d(a_{2},a_{n-2})= 1$,$d(a_{1},a_{n-1})= 1$ these are two terms, \\
$d(a_{3},a_{n-3})= 1$,$d(a_{3},a_{n-2})= 1$,$d(a_{3},a_{n-1})= 1$ these are three terms, so on\\
$d(a_{\lfloor \frac{n}{2}\rfloor-1},a_{\lfloor \frac{n}{2}\rfloor+1})= d(a_{\lfloor \frac{n}{2}\rfloor-1},a_{\lfloor \frac{n}{2}\rfloor+2})=........=d(a_{\lfloor 
	\frac{n}{2}\rfloor-1},a_{{n-1}})=1$ \\these are $(\lfloor \frac{n}{2}\rfloor-1)$ in all\\
similarly $d(a_{\lfloor \frac{n}{2}\rfloor},a_{\lfloor \frac{n}{2}\rfloor+1})=d(a_{\lfloor\frac{n}{2}\rfloor},a_{\lfloor\frac{n}{2}\rfloor+2})=........=d(a_{\lfloor\frac{n}{2}\rfloor},a_{{n-1}})=1$ \\these are also $(\lfloor\frac{n}{2}\rfloor-1)$ in all\\
so on $d(a_{n-2},a_{n-1})= 1$ this is only one term.\\
 \begin{align*}
   W(\Gamma_{E}[\mathbb{Z_\mathbf{m}}])&= \sum\limits_{x,y\in V} d(x,y) .\\
   &= \sum_{dist=2}d(a_{i},a_{j})+\sum_{dist=1}d(a_{i},a_{j})\\
   &= 2\times[(n-3)+(n-5)+.....+(n-{n-1})]+2\times[1+2+....+(\lfloor\frac{n}{2}\rfloor-1)]\\
    &= 2\times[{(n-1)-2}+{(n-1)-4}+.....+{(n-1)-[2{\lfloor\frac{n}{2}\rfloor-1)}]}\\
    &+2\times[1+2+....+(\lfloor\frac{n}{2}\rfloor-1)]\\
     &= 2\times[{(n-1)-2}+{(n-1)-4}+.....+{(n-1)-[2{\lfloor\frac{n}{2}\rfloor-1)}]}\\
     &+2\times[1+2+....+(\lfloor\frac{n}{2}\rfloor-1)]\\
    &= 2\times[{(n-1)(\lfloor\frac{n}{2}\rfloor-1)}]-[2+4+..+2(\lfloor\frac{n}{2}\rfloor-1)]\\
    &+2\times[1+2+....+(\lfloor\frac{n}{2}\rfloor-1)]\\
        &= 2\times[{(n-1)(\lfloor\frac{n}{2}\rfloor-1)}]-2[1+2+..+(\lfloor\frac{n}{2}\rfloor-1)]\\
        &+2\times[1+2+....+(\lfloor\frac{n}{2}\rfloor-1)]\\
          &= 2[{(n-1)(\lfloor\frac{n}{2}\rfloor-1)}]-4[1+2+..+(\lfloor\frac{n}{2}\rfloor-1)]+2\times[1+2+....+(\lfloor\frac{n}{2}\rfloor-1)]\\
           &= 2[{(n-1)(\lfloor\frac{n}{2}\rfloor-1)}]-2[1+2+..+(\lfloor\frac{n}{2}\rfloor-1)]\\
      &=2\times[{(n-1)(\lfloor\frac{n}{2}\rfloor-1)}]-2\times\biggl[\frac{\lfloor\frac{n}{2}\rfloor([\lfloor\frac{n}{2}\rfloor-1])}{2}\bigg].\\
      &=2\times[{(n-1)(\lfloor\frac{n}{2}\rfloor-1)}]-\biggl[\lfloor\frac{n}{2}\rfloor]([\lfloor\frac{n}{2}\rfloor-1])\bigg].\\
    \therefore W(\Gamma_{E}[\mathbb{Z_\mathbf{m}}]) &=2\times[{(n-1)(\lfloor\frac{n}{2}\rfloor-1)}]-\biggl[\lfloor\frac{n}{2}\rfloor([\lfloor\frac{n}{2}\rfloor-1])\bigg].\\       
 \end{align*}
 \begin{align*}\lfloor\frac{n}{2}\rfloor =\left\{
 \begin{array}{ll}
 \frac{n-1}{2}, & \hbox{${\text{if}},n is odd$;} \\
  \frac{n}{2}, & \hbox{${\text{if}} , n is even$.}
 \end{array} \right.
 \end{align*}
 case{1}:- If n is even\\
 \begin{align*}
 \therefore W(\Gamma_{E}[\mathbb{Z_\mathbf{m}}]) &=2\times[{(n-1)([\frac{n}{2}]-1)}]-\biggl[[\frac{n}{2}]([\frac{n}{2}]-1])\bigg].\\   
     &=2\times[{(n-1)(\frac{n}{2}-1)}]-\biggl[\frac{n}{2}(\frac{n}{2}-1])\bigg].\\
      &=2\times[{(n-1)(\frac{n-2}{2})}]-\biggl[\frac{n}{2}(\frac{n-2}{2})\bigg].\\
      &=(n^2-3n+2)-(\frac{n^2-2n}{4}).\\
       &=(\frac{n-2}{4})[4n-4-n].\\
        &=\frac{(n-2)(3n-4)}{4}.\\
      \end{align*}
  case{2}:- If n is odd\\
 \begin{align*}
 \therefore W(\Gamma_{E}[\mathbb{Z_\mathbf{m}}]) &=2\times[{(n-1)([\frac{n}{2}]-1)}]-\biggl[[\frac{n}{2}]([\frac{n}{2}]-1])\bigg].\\   
 &=2\times[{(n-1)(\frac{n-1}{2}-1)}]-\biggl[\frac{n-1}{2}(\frac{n-1}{2}-1])\bigg].\\
 &=2\times[{(n-1)(\frac{n-1-2}{2})}]-\biggl[\frac{(n-1)(n-1-2)}{2}\bigg].\\
 &=(n-3)[(n-1)-(\frac{n-1}{4})].\\
 &=(\frac{n-3}{4})[4n-4-n+1].\\
 &=\frac{3(n-1)(n-3)}{4}+[\frac{n}{2}](\text{an additional term if n is odd}).\\
  &=\frac{3(n-1)(n-3)}{4}+\frac{n-1}{2}.\\
   &=\frac{(n-1)(3n-7)}{4}
 \end{align*}
 \end{proof}
\end{theorem}

\end{document}